\documentclass{article}
\usepackage[T2A]{fontenc}
\usepackage[cp1251]{inputenc}
\usepackage[english]{babel}
\usepackage{color}
\usepackage{amssymb,amsmath,amsthm,amsfonts}
\begin{document}

\addcontentsline{toc}{section}{Abdullaev R.Z, B.A.~Madaminov
Isometries of generalized $log$-spaces.}

\begin{center}
\bf Isometries of generalized $log$-spaces.

R.Z.~Abdullaev, B.A.~Madaminov

arustambay@yandex.com;  \ aabekzod@mail.ru
\end{center}


\begin{abstract}
In this paper studied isometries of $F$-spaces of  integrable functions with logarithm. In particular, using passports of Boolean algebra, a necessary and sufficient condition of isometry $F$-spaces of integrable functions of logarithm with respect to strictly positive $\sigma$-finite measures is proved.\\
In this work,external, internal and generalized $ log $ algebras are considered separately.
\end{abstract}

\textbf{Кeywords:} Functional spaces; $F$-spaces; $log$-algebras; boolean algebras; complete boolean algebras; homogeneous Boolean algebras; passport of boolean algebra; strictly positive; $\sigma$-finite measures; internal $log$-algebras; external $log$-algebras; generalized $log$-algebras; isomorphisms; isometries.

\section*{Introduction}

Many scientific studies with practical applications are reduced to the problems of studying the integration of measurable functions. These include spaces of functions integrable with $p$-th power  ($L_p$-spaces), Orlicz, Lorentz, Marcinkiewicz spaces, Arens algebras and others.  On the other hand, they use an abstract mathematical apparatus based on the general theory of the measure and the Lebesgue integral, the theory of Banach spaces. In modern mathematics, for all these spaces, it is important to establish the conditions of their isometricity and to describe isometries. For all above  mentioned spaces, this problem is solved.

One of the important classes of Banach function spaces are the spaces $L_{p}(\Omega,\mathcal A, \mu),$ $1\leq p<\infty$ of all functions defined on a measurable space $(\Omega,\mathcal A, \mu),$ whose  $p$-th powers are integrable with respect to a finite or $\sigma$-finite measure $\mu$  (functions coinciding almost everywhere are identified). The study of isometries of $L_p$-spaces was started by Banach in \cite{Bana}, where he described all isometries of the spaces $L_p[0,1],\ \ p\neq2.$

	In \cite{Lam}, J.~Lamperti gave a characterization of all linear isometries for $L_p$-spaces $L_{p}(\Omega,\mathcal A, \mu)$, where $(\Omega,\mathcal A, \mu)$ is an arbitrary space with a finite measure $\mu.$  The final result in this formulation belongs to Yedan\cite{Yea}, who gave a complete description of all isometries between $L_p$-spaces associated with various measures. One of the corollaries of such descriptions of isometries in the spaces $L_{p}(\Omega,\mathcal A, \mu)$ is the establishment of isometry for $L_p$-spaces $L_{p}(\Omega,\mathcal A, \mu)$ and $L_{p}(\Omega,\mathcal A, \nu)$ in the case when the measures $\mu$ and $\nu$ are strictly positive finite measures.

In \cite{Dyke}, $F$-spaces of logarithm-integrable functions $L_{log}$ were introduced, which are analogs of $L_p$-spaces.

Interest in these spaces is due to the fact that they are close to the class of Nevanlinn functions holomorphic in the circle and integrable with the logarithm on the boundary of the circle \cite{Dur} , i.e. satisfying the condition
$$L(f):=sup\frac{1}{2\pi}\int\limits_\Omega log(1+|f(re^{i\theta})|)d\theta<\infty.$$

In this paper, we introduce generalizations of these $log$-spaces. They are called generalized $log-$ spaces. For these spaces, a necessary and sufficient condition for their isometry is established for the case of $\sigma$-finite measures.

Let us give some information from \cite{Vlad}.

\emph{Boolean algebra} is a distributive structure with zero and unit, that are unequal to each other.

Boolean algebra is said to be \emph{complete} if every set of its elements has upper and lower bounds.

A mapping $\varphi$ of a Boolean algebra $\nabla_\mu$ to a Boolean algebra $\nabla_\nu$ is said to be \emph{isomorphism} if it is one-to-one and preserves order.

We say that the pairs $(\nabla_\mu,\mu)$ and $(\nabla_\nu,\nu)$ are \emph{isomorphic} if there is a measure-preserving isomorphism between these algebras, i.e. $\mu(x)=\nu(\alpha(x))$ where $\alpha$  is an isomorphism between the Boolean algebras $\nabla_\mu$ and $\nabla_\nu$.

A measure $\mu$ on a Boolean algebra $\nabla_\mu$ is called \emph{strictly positive} if $\mu(x)=0$ implies that $x=0.$

Let $(\Omega, \mathcal A, \mu)$ be a space with a strictly positive $\sigma$-finite measure $\mu,$ and let $\nabla_\nu$ be the complete Boolean algebra of equivalence classes $e=[A]$ from $\mathcal A$ of all almost everywhere equal sets. It is known that $\widehat{\mu}(e)=\mu(A)$, $e \in\nabla_\mu$, is a strictly positive $\sigma$-finite measure on $\nabla_\mu$. The measure $\widehat{\mu}$ will be denoted by $\mu$.

Denote by $L_0(\nabla,\mu)$ the algebra of $\mu$-equivalent classes of complex-valued measurable functions on $(\Omega, \mathcal A, \mu)$.

\textbf{Definition 1.} Let \cite{Dyke} $$L_{log}(\nabla, \mu)=\{f \in L_0(\nabla,\mu):log(1+|f|)d\mu<\infty.\}$$ $ L_{log}(\nabla, \mu)$ be an $F$-space with respect to the $F$-norm $\|f\|_{log}=\int\limits_\Omega log(1+|f|)d\mu$. In \cite{Dyke}, it is also shown that $ L_{log}(\nabla, \mu)$ is an algebra. We call this algebra the \emph{external} $log$-algebra.

\section*{Isometries of external $log$-algebras constructed with respect to $\sigma$-finite measures}

Let $\nabla$ be an arbitrary complete infinite continuous Boolean algebra, $e \in \nabla$, $\nabla_e=[0,e]=\{g \in \nabla:g \leq e\}$. We denote by $\tau(\nabla_e)$ the minimal cardinality of a set dense in $\nabla_e$ in the $(o)$-topology. An infinite complete Boolean algebra $\nabla$ is said to be homogeneous if $\tau(\nabla_e)=\tau(\nabla_g)$ for any nonzero $e, g \in \nabla.$

Let $\nabla$ be a complete non-atomic Boolean algebra and $\mu$ be a strictly positive countably additive $\sigma$-finite measure on $\nabla$. Then the decomposition of $\nabla$ into homogeneous components is at most countable.

\textbf{Definition 2.} We denote by $\{\nabla_{s_i}\}$ homogeneous components of the Boolean algebra $\nabla$ for which
$$
\tau_{s_i}=\tau(\nabla_{s_i})<\tau_{s_{i+1}}, \ \mu(s_i)=\infty,
$$
and we denote by $\{\nabla_{u_i}\}$ homogeneous components of the Boolean algebra $\nabla$ for which
$$
\tau_{u_i}=\tau(\nabla_{u_i})<\tau_{u_{i+1}}, \ \mu(u_i)<\infty.
$$
Here $s_i$ and $u_i$ are units of the Boolean algebras $\nabla_{u_i}$ and $\nabla_{s_i},$ respectively.

Then the matrix
$$\begin{pmatrix}
  \tau_{s_1} & \tau_{s_2} & \ldots \\
  \tau_{u_1} & \tau_{u_2} & \ldots \\
  \mu_1 & \mu_2 & \ldots
\end{pmatrix}
$$
is uniquely defined, which we call the \emph{passport} of the Boolean algebra $\nabla_\mu$ with $\sigma$-finite measure $\mu.$

In the case of a finite measure, we obtain the definition of the passport of a normed Boolean algebra introduced in Theorem 5, p. 270 \cite{Vlad}.

Let $(G, \mathcal B, \nu)$ be a space with a strictly positive $\sigma$-finite measure $\nu,$ and let $\nabla_\nu$ be the complete Boolean algebra of equivalence classes $f=[B]$ from $\mathcal B$ of all sets equal almost everywhere.

Let $\nabla_\mu$ and $\nabla_\nu$ be complete normed (generally speaking, different) Boolean algebras with strictly positive $\sigma$-finite measures $\mu$ and $\nu,$ respectively.

\textbf{Lemma 1.}  \emph{Let $\nabla_\mu$ and $\nabla_\nu$ are homogeneous and do not coincide, $\mu(\Omega)=\nu(G)=\infty.$ For the pairs $(\nabla_\mu,\mu)$ and $(\nabla_\nu,\nu)$ to be isomorphic it is necessary and sufficient that $\tau(\nabla_{\mu})=\tau(\nabla_{\nu })$.}

{\sf Proof.} The isomorphism of the pairs $(\nabla_\mu,\mu)$ and $(\nabla_\nu,\nu)$ implies the isomorphism of the Boolean algebras $\nabla_\mu$ and $\nabla_\nu.$ It follows from here and Theorem 5, p. 270 \cite{Vlad} that the weights of the homogeneous Boolean algebras $\nabla_\mu$ and $\nabla_\nu$ coincide, i.e. $\tau(\nabla_{\mu})=\tau(\nabla_{\nu}).$

Conversely, let $\tau(\nabla_{\mu})=\tau(\nabla_{\nu}).$ Then there exist such countable disjunct decompositions of homogeneous Boolean algebras $\nabla_\mu$ and $\nabla_\nu$ into components $\nabla_\mu e_i$ and $\nabla_\nu f_i,$ respectively, that $\mu(e_i)=1$ and $\mu(f_i)=1$. Therefore, it follows from the corollary of Theorem 5, p. 271\cite{Vlad} that for any $i\in,$  there is a measure-preserving isomorphism $\alpha_i$ from $\nabla_\mu e_i$ to $\nabla_\nu f_i.$ Then the gluing $\bigoplus\alpha_i$ is a measure-preserving isomorphism from $\nabla_\mu$ to $\nabla_\nu.$
Lemma is proved.\hfill $\Box$

Using this Theorem 1, we obtain its generalization for inhomogeneous Boolean algebras.

\textbf{Lemma 2.}  \emph{Let $\nabla_\mu$ и $\nabla_\nu$ be inhomogeneous and do not coincide, the measures  $\mu$ and $\nu$ take the value infinity on each homogeneous component. The pairs $(\nabla_\mu,\mu)$ and $(\nabla_\nu,\nu)$ are isomorphic if and only if the first rows of Boolean algebra passports $(\nabla_\mu,\mu)$ and $(\nabla_\nu,\nu)$}

\begin{equation*}\left(
                  \begin{array}{ccc}
                    \tau_{s_1} & \tau_{s_2} & \dots \\
                     \tau_{u_1} & \tau_{u_2} & \dots \\
                    \mu_1 & \mu_2 & \dots \\
                  \end{array}
                \right) \ \ \mbox{and} \ \
                \left(
                  \begin{array}{ccc}
                    \tau_{s^{'}_1} & \tau_{s^{'}_2} & \dots \\
                      \tau_{u^{'}_1} & \tau_{u^{'}_2} & \dots \\
                    \nu_1 & \nu_2 & \dots \\
                  \end{array}
                \right)
\end{equation*}
\emph{coincide.}

{\sf Proof.} The condition of Theorem implies that the second and third rows of the passports , что are zero. Therefore, it remains to consider the first rows of passports. By virtue of Theorem 1, Theorem 2 is true on every homogeneous component of the Boolean algebra $\nabla_\mu$. Since an isomorphism transforms homogeneous components into homogeneous ones with equal weights, gluing together the isomorphisms from Theorem 1, we obtain a measure-preserving isomorphism of $(\nabla_\mu,\mu)$ and $(\nabla_\nu,\nu).$
Q.E.D.\hfill $\Box$

From Lemma 2 and Theorem 6, p.273 \cite{Vlad}, we obtain an analogue of these theorems for arbitrary $\sigma$-finite measures.

\textbf{Corollary 1.}   \emph{Let $\nabla_\mu$ и $\nabla_\nu$ be inhomogeneous and do not coincide, the measures  $\mu$ and  $\nu$ are $\sigma$-finite. The pairs $(\nabla_\mu,\mu)$ and $(\nabla_\nu,\nu)$ are isomorphic if and only if the first rows of Boolean algebra passports $(\nabla_\mu,\mu)$ and $(\nabla_\nu,\nu)$}
\begin{equation*}\left(
                  \begin{array}{ccc}
                    \tau_{s_1} & \tau_{s_2} & \dots \\
                     \tau_{u_1} & \tau_{u_2} & \dots \\
                    \mu_1 & \mu_2 & \dots \\
                  \end{array}
                \right) \ \ \mbox{and} \ \
                \left(
                  \begin{array}{ccc}
                    \tau_{s^{'}_1} & \tau_{s^{'}_2} & \dots \\
                      \tau_{u^{'}_1} & \tau_{u^{'}_2} & \dots \\
                    \nu_1 & \nu_2 & \dots \\
                  \end{array}
                \right)
\end{equation*}
\emph{coincide.}

The following theorem gives a necessary and sufficient condition for the *-isomorphism of the algebras $L_{log}(\nabla_\mu)$ and $L_{log}(\nabla_\nu)$ given with the help of $\sigma$-finite measures.

\textbf{Theorem 1.}  \cite{Abdul} \emph{Let $\mu$ and $\nu$ be strictly positive $\sigma$-finite measures on non-atomic complete Boolean algebras $\nabla_\mu$ and $\nabla_\nu,$ respectively. Let
\begin{equation*}\left(
                  \begin{array}{ccc}
                    \tau_{s_1} & \tau_{s_2} & \dots \\
                     \tau_{u_1} & \tau_{u_2} & \dots \\
                    \mu_1 & \mu_2 & \dots \\
                  \end{array}
                \right) \ \ \mbox{и} \ \
                \left(
                  \begin{array}{ccc}
                    \tau_{s^{'}_1} & \tau_{s^{'}_2} & \dots \\
                      \tau_{u^{'}_1} & \tau_{u^{'}_2} & \dots \\
                    \nu_1 & \nu_2 & \dots \\
                  \end{array}
                \right)
\end{equation*}
be the passports of the Boolean algebras $(\nabla_\mu,\mu)$ and $(\nabla_\nu,\nu),$ respectively. Then the following conditions are equivalent:}

$(i)$ *-algebras $L_{log}(\nabla_\mu)$ and $L_{log}(\nabla_\nu)$ are *-isomorphic;

$(ii)$ The first and second rows of the passports $\nabla_{\mu}$ и $\nabla_{\nu}$ coincide, and the sequences  $\frac{\mu_i}{\nu_i}$ and  $\frac{\nu_i}{\mu_i}$ are bounded.

The following theorem gives a necessary and sufficient condition for the $F$-spaces $L_{log}(\nabla_\mu)$ and $L_{log}(\nabla_\nu)$ given by finite measures to be isometric.

\textbf{Lemma 3.} \emph{Let $\nabla_{\mu}$ and $\nabla_{\nu}$ be different complete non-atomic homogeneous Boolean algebras with finite strictly positive measures $\mu$ and $\nu,$ respectively. Then the following conditions are equivalent:}

$(i)$ $L_{log}(\nabla_\mu)$ and $L_{log}(\nabla_\nu)$ are isometric;

$(ii)$ a) Wights $\tau(\nabla_\mu)$ and $\tau(\nabla_\nu)$ coincide;

b)  $\mu(u)=\nu(u^{'})$, where $u$ and $u^{'}$ are units of the Boolean algebras $\nabla_{\mu}$ and $\nabla_{\nu},$ respectively.

{\sf Proof.} $(i)\Rightarrow(ii)$  Let $(ii)a)$ is not satisfied, then the homogeneous Boolean algebras $\nabla_{\mu}$ and $\nabla_{\nu}$ have different weights. Therefore, there is no one-to-one correspondence between them. Hence $F$-spaces $L_{log}(\nabla_\mu)$ and $L_{log}(\nabla_\nu)$ are not isometric.

If condition $(ii)b)$ is not satisfied, then we obtain from Theorem 5 \cite{Abdul} that the Boolean algebras $\nabla_{\mu}$ and $\nabla_{\nu}$ are not isometric.

$(ii)\Rightarrow(i)$ Let condition $(ii)$ be satisfied. Then it follows that there exists a measure-preserving isomorphism from $(\nabla_{\mu}, \mu)$ and $(\nabla_{\nu}, \nu)$ , i.e., $\mu(x)=\nu(\alpha(x))$ for any $x\in\nabla_{\mu} $ \cite{Vlad} (Theorem 5, p. 273). Denote by $J_\alpha$ an isomorphism of the algebra $L_\circ(\nabla_\mu)$ onto $L_\circ(\nabla_\nu)$ such that $J_\alpha(x)=\alpha(x)$ for any $x\in\nabla_ {\mu}.$ We get from (\cite{Abdul}, Proposition 3) that
$$
\|f\|_{log,\mu}= \int\limits_{\Omega} log(1+|f(\omega)|)d\mu
=\int\limits_{\Omega} J_{\alpha}(log(1+|f(\omega)|))d\nu=
$$ $$
=\int\limits_{\Omega} (log(1+|J_{\alpha} f(\omega)|))d\nu=\|J_{\alpha}
f(\omega)\|_{log,\nu}
$$
for any $f\in L_{log}(\nabla_\mu).$

We obtain from here that $J_\alpha$ is a bijective linear isometry from $L_{log}(\nabla_\mu)$ to $L_{log}(\nabla_\nu).$

Lemma is proved. \hfill $\Box$

\textbf{Lemma 4.} \emph{Let $\nabla_{\mu}\neq\nabla_{\nu}$ are complete non-atomic homogeneous Boolean algebras with $\sigma$-finite but not finite strictly positive measures $\mu$ and $\nu$. Then the following conditions are equivalent:}

$(i)$ $L_{log}(\nabla_\mu)$ and $L_{log}(\nabla_\nu)$ are isometric;

$(ii)$ a) The weights $\tau(\nabla_\mu)$ and $\tau(\nabla_\nu)$ coincide.

{\sf Proof.}

$(i)\Longrightarrow(ii).$ Let $L_{log}(\nabla_\mu)$ and $L_{log}(\nabla_\nu)$ be isometric. Then there is a one-to-one correspondence between the Boolean algebras $\nabla_{\mu}$ and $\nabla_{\nu}$, i.e. the weights of these algebras coincide.

$(ii)\Longrightarrow(i).$ Let $\tau(\nabla_\mu)=\tau(\nabla_\nu)$ and $\{\nabla_\mu e_i\}_{n=1}^{\infty}$, $\{\nabla_\nu f_i\}_{n=1}^{\infty}$ are disjunct decompositions of Boolean algebras $\nabla_{\mu}$ and $\nabla_{\nu},$ respectively, while $e_i$ and $f_i$ can be chosen so that $\mu(e_i)=\nu(f_i)=1,$ $i=1,2,...$  For any $i,$ there exists a measure-preserving isomorphism $\varphi_i$ of the Boolean algebra $\nabla_\mu e_i$ onto $\nabla_\nu f_i.$  It is clear that the mapping $\varphi,$ defined by the equality $\varphi(x)=\sum^\infty _{i=1}\varphi_i(x\wedge e_n),$ is a measure-preserving isomorphism from $\nabla_{\mu}$ onto $\nabla_{\nu}.$ As well as in the proof of Lemma 3 $(ii)$ we get from here that $L_{log}(\nabla_\mu)$ and $L_{log}(\nabla_\nu)$ are isometric.

Lemma is proved. \hfill $\Box$

The following theorem gives a complete answer to the question of isometricity of external $log$-algebras. It follows from Lemmas 3 and 4, since any isomorphism, and hence isometry, translates homogeneous components into homogeneous ones.

\textbf{Theorem 2.}  \emph{Let $\nabla_\mu$ and $\nabla_\nu$  be complete inhomogeneous Boolean algebras with $\sigma$-finite strictly positive measures  $\mu$ and $\nu,$ respectively. Then $L_{log}(\nabla_\mu)$ and $L_{log}(\nabla_\nu)$ are isometric if and only if the passports of $\nabla_\mu$ and $\nabla_\nu$ coincide.}

We get from Lemma 2 the following

\textbf{Corollary 2.} If the values of $\mu$ and $\nu$ on all homogeneous components are infinite, then the coincidence of the first lines of the passports of the $F$-spaces $L_{log}(\nabla_\mu)$ and $L_{log}(\nabla_\nu )$ is a necessary and sufficient condition for these spaces to be isometric. Moreover, if $\nabla_\mu=\nabla_\nu,$ then $L_{log}(\nabla_\mu)$ and $L_{log}(\nabla_\nu)$ are always isometric.

\section*{Isometries of generalized $log$-algebras constructed with respect to $\sigma$-finite measures}

Let $\mu,$ $\nu$  be strictly positive  $\sigma$-finite measures on a Boolean algebra  $\nabla_\mu$, $\frac{d\nu}{d\mu}=h\geq0$ is the Radon -- Nikodym derivative of $\nu$ with respect to $\mu,$ i.e.

$$\int\limits_{\Omega}fd\nu=\int\limits_{\Omega}hfd\mu,$$
and since $\mu$ and $\nu$ are strictly positive, the support of $h$ is equal to one.
It is clear that then the external $log$-algebras satisfy the equalities

$$
L_{\log}(\nabla_\nu)=\{f \in L_{0}(\nabla):
\int\limits_{\Omega}\log(1+h|f|)d\nu < + \infty\}=\{f \in L_{0}(\nabla):
\int\limits_{\Omega}h\cdot\log(1+|f|)d\mu<+ \infty\}=L^{\nu}_{\log}(\nabla_\mu), \ \
$$

\begin{equation}\label{eq1}
\|f\|_{log,\nu}=\int\limits_{\Omega}\log(1+|f|)d\nu=\int\limits_{\Omega}h(\log(1+|f|))d\mu=\|f\|^{\nu}_{log,\mu}
\end{equation}

Consider now the following analogue of the space of $log$-integrable measurable functions.

\textbf{Definition 3.}

$$L^{(\nu)}_{\log}(\nabla_\mu)=\left\{f \in L_{0}(\nabla):
\int\limits_{\Omega}\log(1+h|f|)d\mu<+ \infty\right\}. $$

The following assertion implies that this space is an $F$-space with respect to the $F$-norm

$$\|f\|^{(\nu)}_{log,\mu}=\int\limits_{\Omega}\log(1+h|f|)d\mu.$$

We call $L^{(\nu)}_{\log}(\nabla_\mu)$  internal $log$-algebras.

\textbf{Proposition 1.} Let $\mu$ and $\nu$ be strictly positive $\sigma$-finite measures on a Boolean algebra $\nabla_\mu.$ Then the function $\|f\|^{(\nu)}_{log,\mu}$ defined on $L^{(\nu)}_{\log}(\nabla_\mu)$ satisfies the following conditions:

$(i)$. $\|f\|^{(\nu)}_{log,\mu}>0$ for all $0 \neq f \in L_{\log}(\nabla_\mu);$

$(ii)$. $\|\alpha f\|^{(\nu)}_{log,\mu}$ for all $f \in
L_{\log}(\nabla_\mu)$ and real numbers  $\alpha$ $|\alpha|\leq
1;$

$(iii)$. $\lim_{\alpha\to 0}\|\alpha f\|^{(\nu)}_{log,\mu}=0$ for all $f \in
L_{\log}(\nabla_\mu);$

$(iv)$. $\|f+g\|^{(\nu)}_{log,\mu}\leq\|f\|^{(\nu)}_{log,\mu}+\|g\|^{(\nu)}_{log,\mu}$ for all $f,
g \in L_{\log}(\nabla_\mu).$

{\sf Proof.}
$(i).$ $\| f\|^{(\nu)}_{log,\mu}=0\Longrightarrow\int\limits_{\Omega}\log(1+h(x)|f(x)|)d\mu=0\Longrightarrow log(1+h(x)|f(x)|)=0$
$\Longrightarrow1+h(x)|f(x)|=1\Longrightarrow h(x)|f(x)|=0.$

$(ii).$ Since the support of $h$ is equal to one, then we get from the last equality that $f=0.$
If $f\in L^{(\nu)}_{\log}(\nabla_\mu),$  then $hf\in L_{\log}(\nabla_\mu).$  Therefore, using the inequality $(ii)$ from Proposition 1, we obtain
$\| \alpha f\|^{(\nu)}_{log,\mu}=\|\alpha h f\|_{\log,\mu}\leq \|h f\|_{\log,\mu}=\|f\|_{\log,\mu}^{(\nu)}.$

$(iii).$ Using $(iii)$ from Proposition 1, we obtain

$\lim\limits_{\alpha\rightarrow0}\| \alpha
f\|^{(\nu)}_{log,\mu}=\lim\limits_{\alpha\rightarrow0}\int\limits_{\Omega}\log(1+\alpha
h|f|)d\mu=\lim\limits_{\alpha\rightarrow0}\| \alpha
f\|_{log,\mu}=0.$

$(iv).\| f+g\|^{(\nu)}_{log,\mu}=\int\limits_{\Omega}\log(1+h(x)|f(x)+g(x)|)d\mu\leq\int\limits_{\Omega}(\log(1+h(x)|f(x)|)$$
$$+log(1+h(x)|g(x)|))d\mu=\|hf\|_{log,\mu}+\|hg\|_{log,\mu}=\|f\|^{(\nu)}_{log,\mu}+\|g\|^{(\nu)}_{log,\mu}.$

Proposition 1 is proved.\hfill $\Box$

\textbf{Theorem 3.}  \emph{For any strictly positive $\sigma$-finite measures $\mu$ and $\nu,$ the spaces $L_{\log}(\nabla_\mu)$ and $L^{(\nu)}_{\log}(\nabla_\mu)$  are isometric.}

{\sf Proof.}

Indeed, for the measures $\mu$ and $\nu$ the mapping $U:L_{\log}(\nabla_\mu)\longrightarrow L^{(\nu)}_{\log}(\nabla_\mu),$ defined by the equality $U(f)=h^{-1}f,$ $f\in L_{log}(\nabla_\mu),$ is a linear surjective isometry from $L_{\log}(\nabla_\mu) $ and $L^{(\nu)}_{\log}(\nabla_\mu).$

Theorem is proved.\hfill $\Box$

\textsc{Remark.} As it can seen from Lemma 1, $F$-spaces $L_{\log}(\nabla_\mu)$ and $L_{\log}(\nabla_\nu)=L^{(\nu)}_{\log} (\nabla_\mu)$ are generally not isometric.

Let  $\mu,\ \nu_1,\ \nu_2\in M(\nabla),\ \mu(\Omega)=1$,  $\nu_1$ and $\nu_2$ are $\sigma$-finite. Then (1) implies that the equalities

\begin{equation}\label{eq2}
\int\limits_{\Omega}fd\nu_1=\int\limits_{\Omega}h_1fd\mu,
\end{equation}

\begin{equation}\label{eq3}
\int\limits_{\Omega}fd\nu_2=\int\limits_{\Omega}h_2fd\mu
\end{equation}
hold.

\textbf{Definition 4.}  We call the spaces
$$L^{\nu_1(\nu_2)}_{\log}(\nabla_\mu)=\{f \in L_{0}(\nabla):
\int\limits_{\Omega}h_1\log(1+h_2|f|)d\mu<+ \infty\}$$
generalized  $log$-algebras.

It will be shown below that the function $\|f\|^{\nu_1(\nu_2)}_{log,\nu}= \int\limits_{\Omega}h_1\log(1+h_2|f|)d\ mu$ is the $F$-norm on the $F$-space $L^{\nu_1(\nu_2)}_{\log}(\nabla_\mu).$ For the case of $h_1$ we obtain the internal $log$-space $L^{\nu}_{\log}(\nabla_\mu)$, and for the case $h_2,$ we obtain the external $L^{\nu}_{\log}(\nabla_\mu).$

\textbf{Lemma 5.}  \emph{The function $\|f\|^{\nu_1(\nu_2)}_{log,\nu}=\int\limits_{\Omega}h_1\log(1+h_2|f|)d\mu$ is the $F$-norm on the $F$-space $L^{\nu_1(\nu_2)}_{\log}(\nabla_\mu)$.}

{\sf Proof. }

By virtue of \eqref{eq1}, we have

$$
\|f\|^{\nu_1(\nu_2)}_{log,\nu}=\int\limits_{\Omega}h_1\log(1+h_2|f|)d\mu=\int\limits_{\Omega}\log(1+h_2|f|)d\nu_1=\|f\|^{(\nu_2)}_{log,\nu_1}.
$$

By virtue of Proposition 1 [\cite{Dyke}, Lemma 2.1], the function $\|f\|^{(\nu_2)}_{log,\mu}$ is the $F$-norm on $L^{(\nu_2)}_{\log}(\nabla_\mu)$, hence, $\|f\|^{\nu_1(\nu_2)}_{log,\nu_1}$ is also the $F$-norm on $L^{\nu_1(\nu_2)}_{\log}(\nabla_\mu)$. \hfill $\Box$

Let $\nu_1,$ $\nu_2,$ $\nu_3,$ $\nu_4,$ $\mu$ be strictly positive finite $\sigma$-finite measures on $\nabla_\mu,$  $\mu(\Omega)=1.$

\textbf{Corollary 3.}  \emph{Let $\nabla_\mu$ be an inhomogeneous Boolean algebra, $\nu_1$ and $\nu_2$ are infinite on every homogeneous component. Then the $F$-spaces $L^{\nu_1(\nu_2)}_{\log}(\nabla_\mu)$ and $L^{\nu_3(\nu_4)}_{\log}(\nabla_\mu)$ are isometric.}

{\sf Proof. }

By virtue of \eqref{eq2},  \eqref{eq3}  we have

$$L^{\nu_1(\nu_2)}_{\log}(\nabla_\mu)=\{f \in L_{0}(\nabla):
\int\limits_{\Omega}h_1\log(1+h_2|f|)d\mu<+ \infty\}=\{f \in L_{0}(\nabla):
\int\limits_{\Omega}\log(1+h_2|f|)d\nu_1<+ \infty\}=L^{(\nu_2)}_{\log}(\nabla_{\nu_1}),$$

$$L^{\nu_3(\nu_4)}_{\log}(\nabla_\mu)=\{f \in L_{0}(\nabla):
\int\limits_{\Omega}h_3\log(1+h_4|f|)d\mu<+ \infty\}=\{f \in L_{0}(\nabla):
\int\limits_{\Omega}\log(1+h_4|f|)d\nu_3<+ \infty\}=L^{(\nu_4)}_{\log}(\nabla_{\nu_3}),$$

$$\|f\|^{\nu_1(\nu_2)}_{log,\nu}=\int\limits_{\Omega}h_1\log(1+h_2|f|)d\mu=\int\limits_{\Omega}\log(1+h_2|f|)d\nu_1=\|f\|^{(\nu_2)}_{log,\nu_1},$$

$$\|f\|^{\nu_3(\nu_4)}_{log,\nu}=\int\limits_{\Omega}h_3\log(1+h_4|f|)d\mu=\int\limits_{\Omega}\log(1+h_4|f|)d\nu_3=\|f\|^{(\nu_4)}_{log,\nu_3}.$$

We obtain from equalities \eqref{eq2},  \eqref{eq3} that $L^{\nu_1(\nu_2)}_{\log}(\nabla_\mu)=L^{(\nu_2)}_{\log}(\nabla_{\nu_1})$ and $L^{\nu_3(\nu_4)}_{\log}(\nabla_\mu)=L^{(\nu_4)}_{\log}(\nabla_{\nu_3}).$  By virtue of Theorem 3, $F$-spaces $L^{(\nu_2)}_{\log}(\nabla_{\nu_1})$ and $L^{(\nu_4)}_{\log}(\nabla_{\nu_3})$ are isometric. Hence, $L^{\nu_1(\nu_2)}_{\log}(\nabla_\mu)$ and $L^{\nu_3(\nu_4)}_{\log}(\nabla_\mu)$ are also isometric since \eqref{eq1} implies that their  $F$-norms coincide.

Corollary is proved.\hfill $\Box$

Combining Theorem 4 \cite{Mad} and Corollary 3, we obtain

\textbf{Theorem 4.}  \emph{Let $\nabla_{\nu_1}=\nabla_{\nu_3}$ be inhomogeneous Boolean algebras. Then the generalized $F$-spaces $L^{\nu_1(\nu_2)}_{\log}(\nabla_\mu)$ and $L^{\nu_3(\nu_4)}_{\log}(\nabla_ \mu)$ are isometric if and only if the third rows of the three-line Boolean algebra passports $\nabla_{\nu_1}$ and $\nabla_{\nu_3}$ coincide.}

R. Abdullaev

Tashkent University of Information Technologies,

Tashkent, 100200, Uzbekistan,

e-mail arustambay@yandex.com

B. Madaminov

Urganch State Pedagogical Institute, Urganch , Uzbekistan,

e-mail aabekzod@mail.ru

\end{document}